\documentclass[12pt,reqno]{amsart}
\usepackage{amsmath, amsfonts, amssymb, amsthm, hyperref}
\usepackage{bm}
\def\l{\left}
\def\r{\right}
\def\bg{\bigg}
\def\({\bg(}
\def\){\bg)}

\def\f{\frac}

\def\eq{\equiv}

\def\Z{\mathbb Z}
\def\C{\mathbb C}
\def\N{\mathbb N}

\def\1{{\bf 1}}

\theoremstyle{plain}
\newtheorem{theorem}{Theorem}[section]

\theoremstyle{lemma}
\newtheorem{lemma}{Lemma}

\newtheorem{conjecture}[theorem]{Conjecture}
\theoremstyle{definition}
\newtheorem*{Ack}{Acknowledgments}

\theoremstyle{Conjecture}

\def\<{\langle}
\def\>{\rangle}

\begin{document}
\hbox{Preprint}
\title{Proof of a congruence concerning truncated hypergeometric series ${}_6F_5$}

\begin{abstract}
In this paper, we mainly prove the following congruence conjectured by J.-C. Liu:
$$
{}_6F_5\bigg[\begin{matrix}\f{5}{4}&\f{1}{2}&\f{1}{2}&\f{1}{2}&\f{1}{2}&\f{1}{2}\\&\f{1}{4}&1&1&1&1\end{matrix}\bigg|\ -1\bigg]_{\f{p-1}{2}}\eq-\f{p^3}{16}\Gamma_p\l(\f{1}{4}\r)^4\pmod{p^5},
$$
where $p\geq5$ are primes with $p\eq3\pmod{4}$.
\end{abstract}
\author{Chen Wang}
\address{Department of Mathematics, Nanjing University, Nanjing 210093,
People's Republic of China}
\email{cwang@smail.nju.edu.cn}

\thanks{
\newline \indent Supported by the National Natural Science Foundation of China (Grant No. 11571162).}
\maketitle

\section{Introduction}
\setcounter{lemma}{0} \setcounter{theorem}{0}
\setcounter{equation}{0}

Define the truncated hypergeometric series
$$
{}_nF_{n-1}\bigg[\begin{matrix}x_1&x_2&\cdots&x_n\\&y_1&\cdots&y_{n-1}\end{matrix}\bigg|\ z\bigg]_m:=\sum_{k=0}^{m}\f{(x_1)_k(x_2)_k\cdots(x_n)_k}{(y_1)_k\cdots(y_{n-1})_k}\cdot\f{z^k}{k!},
$$
where
$$
(x)_k=\begin{cases}x(x+1)\cdots(x+k-1),&\quad k\geq1,\\ 1,&\quad k=0.\end{cases}
$$
Clearly, the truncated hypergeometric series is a finite analogue of the classical hypergeometric series
$$
{}_nF_{n-1}\bigg[\begin{matrix}x_1&x_2&\cdots&x_n\\&y_1&\cdots&y_{n-1}\end{matrix}\bigg|\ z\bigg]:=\sum_{k=0}^{\infty}\f{(x_1)_k(x_2)_k\cdots(x_n)_k}{(y_1)_k\cdots(y_{n-1})_k}\cdot\f{z^k}{k!}.
$$
In 1997, van Hamme \cite{vH} posed several conjectures involving $p$-adic analogue of Ramanujan type series. These congruences are closely related to truncated hypergeometric series. For example, in his paper van Hamme conjectured that
\begin{equation}\label{vanh}
{}_6F_5\bigg[\begin{matrix}\f{5}{4}&\f{1}{2}&\f{1}{2}&\f{1}{2}&\f{1}{2}&\f{1}{2}\\&\f{1}{4}&1&1&1&1\end{matrix}\bigg|\ -1\bigg]_{\f{p-1}{2}}\eq\begin{cases}-p\Gamma_p\l(\f{1}{4}\r)^4\pmod{p^3},&\quad p\eq1\pmod{4},\\ 0\pmod{p^3},&\quad p\eq3\pmod{4},\end{cases}
\end{equation}
where $\Gamma_p(\cdot)$ denotes the $p$-adic Gamma function. Note that the above congruence is a $p$-adic analogue of the following hypergeometric identity due to Ramanujan:
$$
{}_6F_5\bigg[\begin{matrix}\f{5}{4}&\f{1}{2}&\f{1}{2}&\f{1}{2}&\f{1}{2}&\f{1}{2}\\&\f{1}{4}&1&1&1&1\end{matrix}\bigg|\ -1\bigg]=\f{2}{\Gamma\l(\f{3}{4}\r)^4},
$$
which was confirmed by Hardy \cite{H} and Watson \cite{W}. The conjectural congruence of van Hamme was later proved by McCarthy and Osburn \cite{MO}. Note that a lot of congruences involving truncated hypergeometric series have been studied during the past years. One can refer to \cite{DFLST,L,LR,MP,S11a,S11b,S,WP} for details.

In 2015, H. Swisher \cite{S} showed that the congruence \eqref{vanh} also holds modulo $p^5$ for primes $p\eq1\pmod{4}$. Recently, J.-C. Liu \cite{L} investigated \eqref{vanh} modulo $p^5$ for $p\eq3\pmod{4}$ and posed the following conjecture.
\begin{conjecture}\cite{L}\label{conj1}
For primes $p\geq5$ with $p\eq3\pmod{4}$, we have
$$
{}_6F_5\bigg[\begin{matrix}\f{5}{4}&\f{1}{2}&\f{1}{2}&\f{1}{2}&\f{1}{2}&\f{1}{2}\\&\f{1}{4}&1&1&1&1\end{matrix}\bigg|\ -1\bigg]_{\f{p-1}{2}}\eq-\f{p^3}{16}\Gamma_p\l(\f{1}{4}\r)^4\pmod{p^5}.
$$
\end{conjecture}
In the same paper, by using the Mathematica package \verb"Sigma" Liu proved that Conjecture \ref{conj1} holds modulo $p^4$.

In this paper we shall give a complete proof of Conjecture \ref{conj1}.
\begin{theorem}\label{theorem1}
Conjecture \ref{conj1} is true.
\end{theorem}
The organization of this paper is as follows. In next section, we shall give some necessary lemmas in order to show Theorem \ref{theorem1}. We give the proof of Theorem \ref{theorem1} in Section 3.

\medskip
\section{Preliminary results}
Recall that for any $z\in\C$ with $\Re z>0$, the well-known Gamma function is defined as
$$
\Gamma(z):=\int_0^{+\infty}t^{z-1}e^{-t}dt.
$$
One may check that
\begin{equation}\label{gamma1}
\Gamma(z+1)=z\Gamma(z).
\end{equation}
For Gamma function we have the following remarkable reflection formula and duplication formula
\begin{gather}
\label{gamma2}\Gamma(z)\Gamma(1-z)=\f{\pi}{\sin(\pi z)},\\
\label{gamma3}\Gamma(z)\Gamma\l(z+\f{1}{2}\r)=2^{1-2z}\sqrt{\pi}\Gamma(2z).
\end{gather}

Now we recall the $p$-adic Gamma function $\Gamma_p$ which was first introduced by Y. Morita \cite{M} in 1975 as a $p$-adic analogue of the classical Gamma function. Suppose that $p$ is an odd prime. Let $\Z_p$ denote the ring of all $p$-adic integers and let $|\cdot|_p$ denote the $p$-adic norm over $\Z_p$. For each integer $n\geq1$, define
$$
\Gamma_p(n):=(-1)^n\prod_{\substack{1\leq k<n\\ (k,p)=1}}k.
$$
In particular, set $\Gamma_p(0)=1$. For any $x\in\Z_p$, define
$$
\Gamma_p(x):=\lim_{\substack{n\in\N\\ |x-n|_p\rightarrow0}}\Gamma_p(n).
$$
For $p$-adic Gamma function we have the following known results:
\begin{equation}\label{padic1}
\f{\Gamma_p(x+1)}{\Gamma_p(x)}=\begin{cases}-x,&\quad x\not\in p\Z_p,\\ -1,&\quad x\in p\Z_p,\end{cases}
\end{equation}
and
\begin{equation}\label{padic2}
\Gamma_p(x)\Gamma_p(1-x)=(-1)^{\<x\>_p},
\end{equation}
where $\<x\>_p$ denotes the least nonnegative residue of $x$ modulo $p$.

Let $G_k(x)=\Gamma_p^{(k)}(x)/\Gamma_p(x)$, where $\Gamma_p^{(k)}$ denotes the $k$th derivative of $\Gamma_p$. By taking Taylor expansion at $x_0$, we have
$$
\Gamma_p(x)=\Gamma_p(x_0)\sum_{k\geq0}\f{G_k(x_0)}{k!}(x-x_0)^k.
$$
For more properties of $p$-adic Gamma functions one may consult \cite{LR,M,MP}.

To show Theorem \ref{theorem1} we also need the following lemmas.
\begin{lemma}\cite[page 147]{AAR}\label{lemma1}
Let $a,b,c,d,e\in\C$. Then we have the following identity.
\begin{align*}
&{}_6F_5\bigg[\begin{matrix}\f{a}{2}+1&a&b&c&d&e\\ &\f{a}{2}&1+a-b&1+a-c&1+a-d&1+a-e\end{matrix}\bigg| -1\bigg]\\
&\quad\quad=\f{\Gamma(1+a-d)\Gamma(1+a-e)}{\Gamma(1+a)\Gamma(1+a-d-e)}\cdot {}_3F_2\bigg[\begin{matrix}1+a-b-c&d&e\\ &1+a-b&1+a-c\end{matrix}\bigg|\ 1\bigg].
\end{align*}
\end{lemma}

We also need the following identity.
\begin{lemma}\cite[page 149]{AAR}\label{lemma2} For $a,b,c,d,e\in\C$ we have
\begin{align*}
&{}_3F_2\bigg[\begin{matrix}a&b&c\\ &d&e\end{matrix}\bigg|\ 1\bigg]=\f{\pi\Gamma(d)\Gamma(e)}{2^{2c-1}\Gamma(\f{a+d}{2})\Gamma(\f{a+e}{2})\Gamma(\f{b+d}{2})\Gamma(\f{b+e}{2})},
\end{align*}
provided $a+b=1$ and $d+e=2c+1$.
\end{lemma}
\medskip
\section{Proof of Theorem \ref{theorem1}}
\setcounter{lemma}{0} \setcounter{theorem}{0}
\setcounter{equation}{0}
Set
$$
\Psi(x):={}_6F_5\bigg[\begin{matrix}\f{5}{4}&\f{1}{2}&\f{1-ix}{2}&\f{1+ix}{2}&\f{1-x}{2}&\f{1+x}{2}\\&\f{1}{4}&1+\f{ix}{2}&1-\f{ix}{2}&1+\f{x}{2}&1-\f{x}{2}\end{matrix}\bigg| -1\bigg]_{\f{p-1}{2}}.
$$
We first show that
$$
\Psi(p)\eq{}_6F_5\bigg[\begin{matrix}\f{5}{4}&\f{1}{2}&\f{1}{2}&\f{1}{2}&\f{1}{2}&\f{1}{2}\\&\f{1}{4}&1&1&1&1\end{matrix}\bigg| -1\bigg]_{\f{p-1}{2}}\pmod{p^5}.
$$
It is clear that $\Psi(x)$ is a rational function in $x^4$ in $\Z_p[[x^4]]$, so by Taylor expansion we have $\Psi(x)=\sum_{n\geq0}a_nx^{4n}$, where $a_n\in\Z_p$. Define
$$
\Phi(x):={}_6F_5\bigg[\begin{matrix}\f{5-p}{4}&\f{1-p}{2}&\f{1-ix}{2}&\f{1+ix}{2}&\f{1-x}{2}&\f{1+x}{2}\\&\f{1-p}{4}&1+\f{ix-p}{2}&1-\f{ix+p}{2}&1+\f{x-p}{2}&1-\f{x+p}{2}\end{matrix}\bigg| -1\bigg].
$$
Obviously, $\Phi(x)$ is also a rational function in $x^4$ in $\Z_p[[x^4]]$. It is easy to see that $\Phi(x)$ and $\Psi(x)$ share the same coefficients in $\Z_p[[x^4]]$ modulo $p$. On the other hand, by Lemma \ref{lemma1} we have
$$
\Phi(x)=\f{\Gamma(1+\f{x-p}{2})\Gamma(1-\f{x+p}{2})}{\Gamma(\f{3-p}{2})\Gamma(\f{1-p}{2})}\cdot{}_3F_2\bigg[\begin{matrix}\f{1-p}{2}&\f{1-x}{2}&\f{1+x}{2}\\&1+\f{ix-p}{2}&1-\f{ix+p}{2}\end{matrix}\bigg|\ 1\bigg].
$$
Since $(3-p)/2$ and $(1-p)/2$ are all negative integers, $\Phi(x)=0$, that is, its all coefficients vanish in $\Z_p[[x^4]]$. Thus we have $a_1\eq0\pmod{p}$. This concludes that
\begin{equation}\label{key}
\Psi(p)\eq a_0=\Psi(0)={}_6F_5\bigg[\begin{matrix}\f{5}{4}&\f{1}{2}&\f{1}{2}&\f{1}{2}&\f{1}{2}&\f{1}{2}\\&\f{1}{4}&1&1&1&1\end{matrix}\bigg| -1\bigg]_{\f{p-1}{2}}\pmod{p^5}.
\end{equation}

Also, by Lemma \ref{lemma1} and Lemma \ref{lemma2} we have
\begin{align*}
\Psi(p)=&\f{\Gamma\l(1+\f{p}{2}\r)\Gamma\l(1-\f{p}{2}\r)}{\Gamma\l(\f{3}{2}\r)\Gamma\l(\f{1}{2}\r)}\cdot{}_3F_2\bigg[\begin{matrix}\f{1}{2}&\f{1-p}{2}&\f{1+p}{2}\\ &1+\f{ip}{2}&1-\f{ip}{2}\end{matrix}\bigg|\ 1\bigg]\\
=&\f{\Gamma\l(1+\f{p}{2}\r)\Gamma\l(1-\f{p}{2}\r)}{\Gamma\l(\f{3}{2}\r)\Gamma\l(\f{1}{2}\r)}\cdot\f{\pi\Gamma\l(1+\f{ip}{2}\r)\Gamma(1-\f{ip}{2})}{\Gamma\l(\f{3-p-ip}{4}\r)\Gamma\l(\f{3-p+ip}{4}\r)\Gamma\l(\f{3+p-ip}{4}\r)\Gamma\l(\f{3+p+ip}{4}\r)}.
\end{align*}
In view of \eqref{gamma3} we have
$$
\Gamma\l(1+\f{ip}{2}\r)\Gamma\l(1-\f{ip}{2}\r)=\f{1}{\pi}\Gamma\l(\f{1}{2}+\f{ip}{4}\r)\Gamma\l(\f{1}{2}-\f{ip}{4}\r)\Gamma\l(1+\f{ip}{4}\r)\Gamma\l(1-\f{ip}{4}\r).
$$
Thus we deduce that
\begin{equation}\label{Psip}
\Psi(p)=\f{\Gamma\l(1+\f{p}{2}\r)\Gamma\l(1-\f{p}{2}\r)\Gamma\l(\f{1}{2}+\f{ip}{4}\r)\Gamma\l(\f{1}{2}-\f{ip}{4}\r)\Gamma\l(1+\f{ip}{4}\r)\Gamma\l(1-\f{ip}{4}\r)}{\Gamma\l(\f{3}{2}\r)\Gamma\l(\f{1}{2}\r)\Gamma\l(\f{3-p-ip}{4}\r)\Gamma\l(\f{3-p+ip}{4}\r)\Gamma\l(\f{3+p-ip}{4}\r)\Gamma\l(\f{3+p+ip}{4}\r)}.
\end{equation}
Now by \eqref{padic1},
\begin{equation}\label{liuconj1}
\f{\Gamma\l(1+\f{p}{2}\r)\Gamma\l(1-\f{p}{2}\r)}{\Gamma\l(\f{3}{2}\r)\Gamma\l(\f{1}{2}\r)}=\f{p}{2}\cdot\f{\Gamma_p\l(1+\f{p}{2}\r)\Gamma_p\l(1-\f{p}{2}\r)}{\Gamma_p\l(\f{3}{2}\r)\Gamma_p\l(\f{1}{2}\r)}.
\end{equation}
Since $p\eq3\pmod{4}$ we have
\begin{align}
\label{liuconj2}\f{\Gamma\l(1+\f{ip}{4}\r)\Gamma\l(1-\f{ip}{4}\r)}{\Gamma\l(\f{3-p+ip}{4}\r)\Gamma\l(\f{3-p-ip}{4}\r)}=&\prod_{k=(3-p)/4}^0\l(k+\f{ip}{4}\r)\l(k-\f{ip}{4}\r)=\prod_{k=(3-p)/4}^0\l(k^2+\f{p^2}{16}\r),\\
\label{liuconj3}\f{\Gamma\l(\f{1}{2}+\f{ip}{4}\r)\Gamma\l(\f{1}{2}-\f{ip}{4}\r)}{\Gamma\l(\f{3+p+ip}{4}\r)\Gamma\l(\f{3+p-ip}{4}\r)}=&1\bigg/\prod_{k=0}^{(p-3)/4}\l(k+\f{1}{2}+\f{ip}{4}\r)\l(k+\f{1}{2}-\f{ip}{4}\r)\notag\\
=&1\bigg/\prod_{k=0}^{(p-3)/4}\l(\bigg(k+\f{1}{2}\bigg)^2+\f{p^2}{16}\r).
\end{align}
Hence substituting \eqref{liuconj1}--\eqref{liuconj3} into \eqref{Psip} and noting that
$$
\Gamma_p\l(1+\f{p}{2}\r)\Gamma_p\l(1-\f{p}{2}\r)\eq\Gamma_p(1)^2\l(1+\f{p}{2}G_1(1)\r)\l(1-\f{p}{2}G_1(1)\r)\eq\Gamma_p(1)^2\pmod{p^2}
$$
and
$$
\Gamma_p\l(\f{p+1}{4}\r)\Gamma_p\l(\f{1-p}{4}\r)\eq\Gamma_p\l(\f{1}{4}\r)^2\l(1+\f{p}{4}G_1\l(\f{1}{4}\r)\r)\l(1-\f{p}{4}G_1\l(\f{1}{4}\r)\r)\eq\Gamma_p\l(\f{1}{4}\r)^2\pmod{p^2},
$$
we have
\begin{align*}
\Psi(p)=&\f{p^3}{32}\cdot\f{\Gamma_p\l(1+\f{p}{2}\r)\Gamma_p\l(1-\f{p}{2}\r)}{\Gamma_p\l(\f{3}{2}\r)\Gamma_p\l(\f{1}{2}\r)}\cdot\f{\prod_{k=(3-p)/4}^{-1}(k^2+p^2/16)}{\prod_{k=0}^{(p-3)/4}((k+1/2)^2+p^2/16)}\\
\eq&-\f{p^3}{16}\cdot\f{\prod_{k=(3-p)/4}^{-1}k^2}{\prod_{k=0}^{(p-3)/4}(k+1/2)^2}=-\f{p^3}{16}\cdot\f{\Gamma_p((p+1)/4)^2\Gamma_p(1/2)^2}{\Gamma_p((p+3)/4)^2}\\
\eq&-\f{p^3}{16}\cdot\Gamma_p\l(\f{p+1}{4}\r)^2\Gamma_p\l(\f{1-p}{4}\r)^2\eq-\f{p^3}{16}\cdot\Gamma_p\l(\f{1}{4}\r)^4\pmod{p^5}.
\end{align*}
Combining this with \eqref{key} we immediately obtain the desired theorem.\qed

\begin{Ack}
The author would like to thank Prof. Hao Pan at Nanjing University of Finance and Economics and Dr. Hai-Liang Wu at Nanjing University for their helpful comments.
\end{Ack}

\end{document}